\documentclass[a4paper]{amsart}

\usepackage{fullpage}
\usepackage{xcolor}
\usepackage{amsmath,amssymb,amsthm}
\usepackage{fullpage}
\usepackage{graphicx}  
\usepackage[colorlinks=true]{hyperref}
\hypersetup{urlcolor=blue, citecolor=blue}

\theoremstyle{definition}

\newtheorem{remark}{Remark}

\usepackage{graphicx}
\usepackage{stackrel}
\usepackage{tikz} 
\usetikzlibrary{matrix,arrows}
\usepackage{color}
\usepackage{textcomp}
\usepackage{enumerate}
\usepackage{mathtools}
\usepackage{subcaption}
\usepackage{float}
\usepackage{mathrsfs}
\usepackage{multirow}
\usepackage{multicol}  
\usepackage{tabularx}
\usepackage{booktabs}
\usepackage{comment}
\usepackage{xypic}
\usepackage{bbm}
\usepackage[all]{xy}
\usepackage{algorithmic}
\usepackage[ruled,vlined]{algorithm2e}
\usepackage{setspace}

\newcommand{\R}{\ensuremath{\mathbb{R}}}

\newcommand{\cP}{\ensuremath{\mathcal{P}}}

\newcommand{\tr}{\ensuremath{\mathrm{Tr}}}

\DeclareMathOperator*{\argmin}{argmin}
\DeclareMathOperator*{\rank}{rank}
\DeclareMathOperator*{\diag}{diag}

\begin{document}

\title{An Adaptation for Iterative Structured Matrix Completion}

\author{Henry Adams}
\email{henry.adams@colostate.edu}
\author{Lara Kassab}
\email{kassab@math.colostate.edu}
\author{Deanna Needell}
\email{deanna@math.ucla.edu}



\thanks{Needell was partially supported by NSF CAREER $\#1348721$ and NSF BIGDATA $\#1740325$.}

\thanks{Corresponding author: Lara Kassab}

\maketitle

\begin{abstract}
The task of predicting missing entries of a matrix, from a subset of known entries, is known as \textit{matrix completion}. 
In today's data-driven world, data completion is essential whether it is the main goal or a pre-processing step.
Structured matrix completion includes any setting in which data is not missing uniformly at random.
In recent work, a modification to the standard nuclear norm minimization (NNM) for matrix completion has been developed to take into account \emph{sparsity-based} structure in the missing entries.
This notion of structure is motivated in many settings including recommender systems, where the probability that an entry is observed depends on the value of the entry.
We propose adjusting an Iteratively Reweighted Least Squares (IRLS) algorithm for low-rank matrix completion to take into account sparsity-based structure in the missing entries.
We also present an iterative gradient-projection-based implementation of the algorithm that can handle large-scale matrices.
Finally, we present a robust array of numerical experiments on matrices of varying sizes, ranks, and level of structure.
We show that our proposed method is comparable with the adjusted NNM on small-sized matrices, and often outperforms the IRLS algorithm in structured settings on matrices up to size $1000 \times 1000$.
\end{abstract}

\section{Introduction}
\textit{Matrix completion} is the task of filling-in, or predicting, the missing entries of a partially observed matrix from a subset of known entries.
In today's data-driven world, data completion is essential, whether it is the main goal as in recommender systems, or a pre-processing step for other tasks like regression or classification.
One popular example of a data completion task is the \textit{Netflix Problem}~\cite{bell2007lessons,bennett2007netflix,koren2009matrix}, which was an open competition for the best collaborative filtering algorithm to predict unseen user ratings for movies.
Given a subset of  user-movie ratings, the goal is to predict the remaining ratings, which can be used to decide whether a certain movie should be recommended to a user.
The Netflix Problem can be viewed as a matrix completion problem where the rows represent users, the columns represent movies, and the entries of the matrix are the corresponding user-movie ratings, most of which are missing.

Matrix completion problems are generally ill-posed without some additional information, since the missing entries could be assigned arbitrary values.
In many instances, the matrix we wish to recover is known to be low-dimensional in the sense that it is low-rank, or approximately low-rank.
For instance, a data matrix of all user-ratings of movies may be approximately low-rank because it is commonly believed that only a few factors contribute to an individual's tastes or preferences~\cite{candes2009exact}.
Low-rank matrix completion is a special case of the \textit{affine rank minimization problem}, which arises often in machine learning, and is known to be NP-hard~\cite{fazel2002matrix,recht2010guaranteed}.

Standard matrix completion strategies typically assume that there are no \textit{structural differences} between observed and missing entries, which is an unrealistic assumption in many settings.
Recent works~\cite{chen2015completing,molitor2018matrix,robin2019main,schnabel2016recommendations,sportisse2018imputation} address various notions of the problem of structured matrix completion.
General notions of structural difference include any setting in which whether an entry is observed or unobserved does not occur uniformly at random.
For example, the probability that an entry is observed could depend not only on the value of that entry, but also on its location in the matrix. 
For instance, certain rows (or columns) may have substantially more entries than a typical row (or column); this happens in the Netflix Problem for very popular movies or so-called ``super-users''.

In~\cite{molitor2018matrix}, Molitor and Needell propose a modification to the \textit{standard nuclear norm minimization} for matrix completion to take into account structure when the submatrix of unobserved entries is sparse, or when the unobserved entries have lower magnitudes than the observed entries~\cite{molitor2018matrix}.
In our work, we focus on this notion of structure, in which the probability that an entry is observed or not depends mainly on the value of the entry.
In particular, we are interested in sparsity-based structure in the missing entries, whereby the submatrix of missing values is close to 0 in the $L_1$ or $L_0$ norm sense.
This is motivated by many situations in which the missing values tend to be near a certain value.
For instance, missing data in chemical measurements might indicate that the measurement value is lower than the limit of detection of the device, and thus a typical missing measurement is smaller in value than a typical observed measurement.
Similarly, in medical survey data, patients are more likely to respond to questions addressing noticeable symptoms, whereas
a missing response may indicate a lack of symptoms~\cite{molitor2018matrix}.
In the Netflix problem, a missing rating of a movie might indicate the user's lack of interest in that movie, thus suggesting a lower rating than otherwise expected.
More generally, in survey data, incomplete data may be irrelevant or unimportant to the individual, therefore suggesting structure in the missing observations~\cite{molitor2018matrix}.
For an example in the setting of sensor networks, suppose we are given partial information about the distances between sensors, where distance estimates are based on signal strength readings~\cite{candes2009exact}, and 
we would like to impute the missing signal strength readings.
Signals may be missing because of low signal strength, indicating that perhaps these sensors are far from each other (or there are geographic obstacles between them).
Thus, we obtain a partially observed matrix with structured observations---missing entries tend to have lower signal strength.
Sensor networks give a low-rank matrix completion problem, perhaps of rank equal to two if the sensors are
located in a plane, or three if they are located in three-dimensional space~\cite{linial1995geometry,so2007theory}.
Therefore, in these settings, we expect that the missing entries admit a sparsity-based structure in the $L_1$ norm sense.

\subsection{Background and Related Work}\label{ss:related work}

The \textit{Affine Rank Minimization Problem} (ARMP), or the problem of finding the minimum rank matrix in an affine set, is expressed as
\begin{equation}
\begin{aligned}
& \underset{X}{\text{minimize}}
& & \rank (X) \\
& \text{subject to}
& & \mathcal A (X) = b,
\end{aligned}
\label{eq:ARMP}
\end{equation}
where matrix $X \in \R^{m\times n}$ is the optimization variable, $\mathcal A \colon \R^{m\times n} \to \R^q$ is a linear map, and $b \in \R^q$ denotes the measurements.
The affine rank minimization problem arises frequently in applications like system identification and control~\cite{liu2009interior}, collaborative filtering, low-dimensional Euclidean embeddings~\cite{fazel2003log}, sensor networks~\cite{biswas2006semidefinite,schmidt1986multiple,singer2008remark}, quantum state tomography~\cite{gross2011recovering,gross2010quantum}, signal processing, and image processing.

Many algorithms have been proposed for ARMP, e.g.\ reweighted nuclear norm minimization~\cite{mohan2010reweighted}, Singular
Value Thresholding (SVT)~\cite{cai2010singular}, Fixed Point Continuation Algorithm (FPCA)~\cite{ma2011fixed}, Iterative
Hard Thresholding (IHT)~\cite{goldfarb2011convergence}, Optspace~\cite{keshavan2009gradient}, Singular Value Projection (SVP)~\cite{jain2010guaranteed}, Atomic Decomposition for Minimum Rank Approximation (AdMiRA)~\cite{lee2010admira}, 
Alternating Minimization approach~\cite{jain2013low},
and the accelerated proximal gradient algorithm for nuclear
norm regularized linear least squares problems (NNLS)~\cite{toh2010accelerated}, etc.

The low-rank matrix completion problem can be formulated as follows~\cite{candes2010power,recht2011simpler}.
Suppose we are given some matrix $X \in \R^{m \times n}$ with a set $\Omega$ of partially observed entries, of size $|\Omega| \ll m n$.
The goal is to recover the missing elements in $X$.
The low-rank matrix completion problem is a special case of the affine rank minimization problem
where the set of affine constraints restrict certain entries of the matrix $X$ to equal observed values.
In this case, the linear operator $\mathcal A$ is a sampling operator, and the problem can be written as
\begin{equation*}
\begin{aligned}
& \underset{X}{\text{minimize}}
& & \rank (X) \\
& \text{subject to}
& & X_{ij} = M_{ij} , \quad (i, j) \in \Omega,
\end{aligned}
\end{equation*}
where $M$ is the matrix we would like to recover, and where $\Omega$ denotes the set of entries which are revealed.
We define the sampling operator $\cP_\Omega(X)\colon \R^{m \times n} \to \R^{m \times n}$ via
\begin{equation*}
(\cP_\Omega(X))_{ij} = \begin{cases} 
      X_{ij} & (i,j) \in \Omega \\
      0 & (i,j) \not \in \Omega, \\
   \end{cases}
 \end{equation*}
as in~\cite{candes2010power}.
Further, $\Omega^c$ denotes the complement of $\Omega$, i.e., all index pairs $(i, j)$ that are not in $\Omega$.
Thus, $\Omega^c$ corresponds to the collection of missing entries.
The \emph{degrees of freedom ratio} of a partially observed $m \times n$ matrix of rank $r$ is given by $FR = r(m+n - r)/|\Omega|$.
Thus, the larger the degrees of freedom ratio is, the harder it becomes to recover the matrix $M$.

The rank minimization problem~\eqref{eq:ARMP} is NP-hard in general, and therefore we consider its convex relaxation~\cite{candes2010matrix,candes2009exact,candes2010power,fazel2002matrix,recht2010guaranteed},
\begin{equation}
\begin{aligned}
& \underset{X}{\text{minimize}}
& & \|X\|_* \\
& \text{subject to}
& & \mathcal A (X) = b,
\end{aligned}
\label{eq:NNM}
\end{equation}
where $\|\cdot \|_{*}$ denotes the nuclear norm, given by the sum of singular values.

Inspired by the iteratively reweighted least squares (IRLS) algorithm for sparse vector recovery analyzed in~\cite{daubechies2010iteratively}, iteratively reweighted least squares algorithms~\cite{fornasier2011low,kummerle2018harmonic,mohan2012iterative} have been proposed as a computationally efficient method for low-rank matrix recovery
(see Section~\ref{ss:related low-rank matrix recovery}). 
Instead of minimizing the nuclear norm, the algorithms essentially minimize the Frobenius norm of the matrix, subject to affine constraints.
Properly reweighting this norm produces low-rank solutions under suitable assumptions.
In~\cite{mohan2012iterative}, Mohan and Fazel propose a family of Iterative Reweighted Least Squares algorithms for matrix rank minimization, called IRLS-$p$ (for $0\leq p \leq 1$), as a computationally efficient way to improve over the performance of nuclear norm minimization.
In addition, a gradient projection algorithm is presented as an efficient implementation for the algorithm, which exhibits improved recovery when compared to existing algorithms.

Generally, standard matrix completion strategies assume that there are no structural differences between observed and unobserved entries.
However, recent works~\cite{chen2015completing,chen2013low,molitor2018matrix,razenshteyn2016weighted,robin2019main,schnabel2016recommendations,sportisse2018imputation} also address various notions of the problem of structured matrix completion in mathematical, statistical, and machine learning frameworks.
In our work, we are interested in sparsity-based structure.
This notion of structure was proposed in~\cite{molitor2018matrix}, where the standard nuclear norm minimization problem for low-rank matrix completion is modified to take into account sparsity-based structure by regularizing the values of the unobserved entries.
We refer to this algorithm as Structured NNM (see Section~\ref{ss:problem statement}).

\subsection{Contribution}\label{ss:contribution}

We adapt an iterative algorithm for low-rank matrix completion to take into account sparsity-based structure in unobserved entries by adjusting the IRLS-$p$ algorithm studied in~\cite{mohan2012iterative}.
We refer to our algorithm as \emph{Structured IRLS}.
We also present a gradient-projection-based implementation, called \emph{Structured sIRLS} (motivated by sIRLS in~\cite{mohan2012iterative}).
Our code is available at~\cite{StructuredIRLScode}.
The main motivations for our approach, along with its advantages, are as follows:
\begin{itemize}

\item \textbf{Iterative algorithm for structured matrix completion.}
Much work has been put into developing iterative algorithms (SVT~\cite{cai2010singular}, FPCA~\cite{ma2011fixed}, IHT~\cite{goldfarb2011convergence}, IRLS~\cite{fornasier2011low,kummerle2018harmonic,mohan2012iterative}, etc.) for ARMP, rather than solving the nuclear norm convex heuristic (NNM).
We develop the first (to our knowledge) iterative algorithm that addresses the structured low-rank matrix completion problem, for which Structured NNM has been proposed.
Indeed, iterative methods are well-known to offer ease of implementation and reduced computational resources, making our approach attractive in the large-scale settings.

\item \textbf{Comparable performance with Structured NNM.}
Structured NNM adapts nuclear norm minimization and $\ell_1$ norm minimization, which are common heuristics for minimizing rank and inducing sparsity, respectively.
For various structured regimes, we consider small-sized matrices and show that our proposed iterative method is comparable to Structured NNM on ``hard'' matrix completion problems and with ``optimal'' parameter choices for Structured NNM.

    
\item \textbf{Improved IRLS recovery for structured matrices.}
We show that in structured settings, Structured sIRLS often performs better than the sIRLS algorithm, as follows.
We perform numerical experiments that consider $19^2=361$ combinations of different sampling rates of the zero and nonzero entries, in order to demonstrate various levels of sparsity in the missing entries.
Consider for example Figure~\ref{fig:exp1000_rank10}, on matrices of size $1000 \times 1000$ of rank 10, in which our proposed method outperforms standard sIRLS in over 90\% of the structured experiments (and also in many of the unstructured experiments).
    
\item \textbf{Handle hard problems.} 
We consider problems of varying degrees of freedom, and a priori rank knowledge.
We show that Structured sIRLS often outperforms the sIRLS algorithm in structured settings for hard matrix completion problems, i.e.\ where the degrees of freedom ratio is greater than 0.4.

\item \textbf{Handle noisy measurements.} 
We consider matrices with noisy measurements with two different levels of noise.
We show that for small enough noise Structured sIRLS often performs better than sIRLS in structured settings.
As the noise gets larger, both converge to the same performance.

\end{itemize}



\subsection{Organization}
We review related iteratively reweighted least squares algorithms for recovering sparse vectors and low-rank matrices in Section~\ref{sec:IRLS-algorithms}.
In Section~\ref{sec:structured IRLS-algorithms}, we describe the structured matrix completion problem, propose for this problem an iterative algorithm, Structured IRLS, and present preliminary analytic remarks.
Furthermore, we present a computationally efficient implementation, Structured sIRLS.
In Section~\ref{sec: numerical experiments}, we run numerical experiments to showcase the performance of this method, and compare it to the performance of
sIRLS and Structured NNM on various structured settings.

\section{Iteratively Reweighted Least Squares Algorithms}\label{sec:IRLS-algorithms}
In this section, we set notation for the rest of the paper, and we review related algorithms for recovering sparse vectors and low-rank matrices.

\subsection{Notation}\label{ss:notation}
The entries of a matrix $X \in \R^{m \times n}$ are denoted by lowercase letters, i.e.\ $x_{ij}$ is the entry in row $i$ and column $j$ of $X$.
Let $I$ denote the identity matrix and $\mathbbm 1$ the vector of all ones.
The \emph{trace} of a square matrix $X \in \R^{m \times m}$ is the sum of its diagonal entries, and is denoted by $\tr(X) = \sum\limits_{i=1}^m x_{ii}$.  
We denote the adjoint matrix of $X$ by $X^* \in \R^{n \times m}$.
Without loss of generality, we assume $m\leq n$ and we write the singular value decomposition of $X$ as
\[ X = U \Sigma V^*.\]
Here $U \in \R^{m \times m}$ and $V \in \R^{n \times n}$ are unitary matrices, and $\Sigma = \diag(\sigma_1, \cdots, \sigma_m) \in \R^{m \times n}$ is a diagonal matrix, where $\sigma_1 \geq \sigma_2 \geq \cdots \geq \sigma_m \geq 0$ are the \textit{singular values}.
The \emph{rank} of $X \in \R^{m \times n}$, denoted by $\rank(X)$, equals the number of nonzero singular values of $X$.
Further, the \emph{Frobenius norm} of the matrix $X$ is defined by
\[\|X\|_F = \sqrt{\tr(XX^*)} = \left( \sum\limits_{i=1}^m \sum\limits_{j=1}^n x_{ij}^2\right)^{1/2} = \left( \sum\limits_{i=1}^{m} \sigma_i^2 \right)^{1/2} .\]
The \emph{nuclear norm} of a matrix $X$ is defined by $\|X\|_* = \sum\limits_{i=1}^{m} \sigma_i$.
Given a vector $w \in \R^n$ of positive weights, we define the \emph{weighted $\ell_2$ norm} of a vector $z \in \R^n$ as
 \[ \|z\|_{\ell_2(w)}  = \left(\sum\limits_{i=1}^{n} w_i z_i^2\right)^{1/2}. \]
Let $z(X)$ denote the vector of missing entries of a matrix $X$, and let $z^2(X)$ denote the corresponding vector with entries squared, i.e.\ $z^2(X) = z(X) \odot z(X) $ where $\odot$ denotes elementwise multiplication.

\subsection{Sparse Vector Recovery} 
\label{ss:related sparse recovery}
Given a vector $x$, the value $\| x \|_0$ denotes the number of nonzero entries of $x$, and is known as the \emph{$\ell_0$ norm} of $x$.
The sparse vector recovery problem is described as
\begin{equation}
\begin{aligned}
& {\text{minimize}}
& & \|x\|_0,\\
& \text{subject to}
& & Ax = b,
\end{aligned}
\label{eq:sparse recovery}
\end{equation}
where $x \in \R^n$ and $A \in \R^{m \times n}$.
This problem is known to be NP-hard.
A commonly used convex heuristic for this problem is $\ell_1$ minimization~\cite{candes2005decoding,donoho1989uncertainty},
\begin{equation}
\begin{aligned}
& {\text{minimize}}
& & \|x\|_1,\\
& \text{subject to}
& & Ax = b.
\end{aligned}
\label{eq:ell_1 minimization}
\end{equation}
Indeed, many algorithms for solving~\eqref{eq:sparse recovery} and~\eqref{eq:ell_1 minimization} have been proposed.
In~\cite{daubechies2010iteratively}, Daubechies et al.\ propose and analyze a simple and computationally efficient reweighted algorithm for sparse recovery, called the Iterative Reweighted Least Squares algorithm, IRLS-$p$, for any $0 < p \leq 1$.
Its $k$-th iteration is given by
\begin{equation*}
    x^{k+1} = \argmin _x \left\{\sum \limits _{i=1} ^n w_i^k x_i^2~:~Ax =b \right\},
\end{equation*}
where $w^k \in \R^n$
is a weight vector with $w^k_i = ( |x_i^k|^2+ \epsilon_k^2)^{p/2-1}$, and where $\epsilon_k > 0$ is a regularization
parameter added to ensure that $w^k$ is well-defined.
For $p = 1$,~\cite{daubechies2010iteratively} gives a theoretical guarantee for sparse recovery similar to $\ell_1$ norm minimization.

\subsection{Low-rank Matrix Recovery}
\label{ss:related low-rank matrix recovery}
We review two related algorithms~\cite{fornasier2011low,mohan2012iterative} for low-rank matrix recovery that generalize the iteratively reweighted least squares algorithm analyzed in~\cite{daubechies2010iteratively} for sparse recovery.
In general, minimizing the Frobenius norm subject to affine constraints does not lead to low-rank solutions; however, properly reweighting this norm produces low-rank solutions under suitable assumptions~\cite{fornasier2011low,mohan2012iterative}.

In~\cite{fornasier2011low}, Fornasier et al.\ propose a variant of the reweighted least squares algorithm for sparse recovery for nuclear norm minimization (or low-rank matrix recovery), called IRLS-M.
The $k$-th iteration of IRLS-M is given by
\begin{equation}
    X^{k+1} = \argmin \limits _{X} \left \{\| (W^ {k})^{1/2} X\|_F^2 \colon  \cP_{\Omega}(X) = \cP_{\Omega}(M)  \right \}.
    \label{alg: IRLS-M}
\end{equation}
Here $W^ k \in \R^{m \times m}$ is a weight matrix defined as $W^0 = I$, and for $k>0$, $W^k = U^k (\Sigma _{\epsilon _ k}^ k) ^{-1} (U^k)^*$, where $X^k (X ^k)^* = U^k (\Sigma^ k) ^2 (U^ k)^*$ and $\Sigma _{\epsilon _ k} = \diag(\max\{\sigma_j , \epsilon_k\})$.
Indeed, each iteration of~\eqref{alg: IRLS-M}  minimizes a weighted Frobenius norm of the matrix X.
Under the assumption that the linear measurements fulfill a suitable generalization of the \textit{Null Space Property} (NSP), the algorithm is guaranteed to iteratively recover any matrix with an error on the order of the best rank $k$ approximation~\cite{fornasier2011low}.
The algorithm essentially has the same recovery guarantees as nuclear norm minimization.
Though the Null Space Property fails in the matrix completion setup, the authors illustrate numerical experiments which show that the IRLS-M algorithm still works very well in this setting for recovering low-rank matrices.
Further, for the matrix completion problem, the algorithm takes advantage of the Woodbury matrix identity, allowing an expedited solution to the least squares problem required at each iteration~\cite{fornasier2011low}.

In~\cite{mohan2012iterative}, Mohan and Fazel propose a related family of Iterative Reweighted Least Squares algorithms for matrix rank minimization, called IRLS-$p$ (for $0\leq p \leq 1$), as a computationally efficient way to improve over the performance of nuclear norm minimization.
The $k$-th iteration of IRLS-$p$ is given by
\begin{equation}
    X^{k+1} = \argmin \limits _{X} \left \{ \tr (W_p^{k}X^\top X) \colon  \cP_{\Omega}(X) = \cP_{\Omega} (M)  \right \},
    \label{alg: IRLS-p}
\end{equation}
where $W _p ^ k \in \R^{m \times m}$ is a weight matrix defined as $W _p ^ 0 = I$, and for $k>0$, $W _p ^ k = ((X^{k})^{\top}X^k + \gamma ^k I)^{\frac{p}{2}-1}$.
Here $\gamma^k > 0$ is a regularization parameter added to ensure that $W _p ^ k$ is well-defined.
Each iteration of~\eqref{alg: IRLS-p}  minimizes a weighted Frobenius norm of the matrix X, since 
\[ \tr (W_p^{k-1}X^\top X) = \|(W_p^{k-1})^{1/2}X\|_F^2.\]
The algorithms can be viewed as (locally) minimizing certain smooth approximations to the rank function.
When $p = 1$, theoretical guarantees are given similar to those for nuclear norm minimization, i.e., recovery of low-rank matrices under the assumptions that the operator defining the constraints satisfies a specific \textit{Null Space Property}.
Further, for $p < 1$, IRLS-$p$
shows better empirical performance in terms of recovering low-rank matrices than nuclear
norm minimization.
In addition, a gradient projection algorithm, IRLS-GP, is presented as an efficient implementation for IRLS-$p$.
Further, this same paper presents a related family of algorithms sIRLS-$p$ (or short IRLS), which can be seen as a first-order method for locally minimizing a smooth approximation to the rank function.
The results exploit the fact that these
algorithms can be derived from the KKT conditions for minimization problems whose
objectives are suitable smooth approximations to the rank function~\cite{mohan2012iterative}.
We will refer to IRLS-$p$ (resp.\ sIRLS-$p$) studied in~\cite{mohan2012iterative} as IRLS (resp.\ sIRLS).

The algorithms proposed in~\cite{fornasier2011low,mohan2012iterative} differ mainly in their implementations, and in the update rules of the weights and their corresponding regularizers.
In IRLS-M~\cite{fornasier2011low}, the weights are updated as $ W^k  = U^k \diag(\max({\sigma_j^k , \epsilon_k })^{-1})(U^k)^*$, 
and in IRLS-$p$~\cite{mohan2012iterative} they are updated as $ W^k  = U^k \diag(((\sigma_j^k) ^2+ \gamma_k ^2)^{-1/2})(U^k)^*$.
Further, each of the regularization parameters $\epsilon_k$ and $\gamma_k$ are updated differently.
The IRLS-M algorithm makes use of the rank of the matrix (either given or estimated), and thus the choice of parameter $\epsilon_k$ depends on this given or estimated rank.
On the other hand, the IRLS-$p$ algorithm chooses and updates its regularizer $\gamma_k$ based on prior sensitivity experiments.

\begin{table}[htbp]
\centering
{ \renewcommand{\arraystretch}{1.15}
\begin{tabular}{ l p{0.66\textwidth}  l }
            \hline
            \multicolumn{2}{c}{Terminology} \\
            \hline
	        NNM \qquad& Nuclear Norm Minimization\\
	        Structured NNM \qquad& Adjusted NNM for sparsity-based structure in the missing entries, proposed in~\cite{molitor2018matrix}\\
            IRLS-$p$ \qquad& Iterative Reweighted Least Squares algorithms for matrix rank minimization, proposed in~\cite{mohan2012iterative} \\
            sIRLS \qquad& short IRLS-$p$,  proposed in~\cite{mohan2012iterative} \\
            Structured IRLS \qquad& Our proposed algorithm: adjusted IRLS for sparsity-based structure in the missing entries \\
            Structured sIRLS \qquad& Our proposed implementation: adjusted sIRLS for sparsity-based structure in the missing entries\\
            \hline
\end{tabular}
}
\end{table}

\section{Structured Iteratively Reweighted Least Squares Algorithms}\label{sec:structured IRLS-algorithms}

In this section, we first introduce the structured matrix completion problem. 
Second, we introduce and analyze our proposed algorithm and implementation.

\subsection{Problem Statement}
\label{ss:problem statement}
In~\cite{molitor2018matrix}, the authors propose adjusting
the standard nuclear norm minimization (NNM) strategy for matrix
completion to account for structural differences between
observed and unobserved entries.
This could be achieved by adding to problem~\eqref{eq:NNM} a regularization term on the unobserved entries, which results in a semidefinite program:
\begin{equation}
\begin{aligned}
& \underset{X}{\text{minimize}}
& & \|X\|_* + \alpha \|\cP_{\Omega^c} (X)\| \\
& \text{subject to}
& & \cP_\Omega (X) = \cP_\Omega (M),
\end{aligned}
\label{eq:structured nnm}
\end{equation}
where $\alpha >0$, and where $\| \cdot \|$ is an appropriate matrix norm.
If most of the missing entries are zero except for a few, then the $\ell_1$ norm is a natural choice\footnote{The method can be rescaled if there instead is a preference for the missing entries to be near a nonzero constant.
}.
If the missing entries are mostly close to zero, then the $\ell_2$ norm is a natural choice.
The authors show that the proposed method outperforms nuclear norm minimization in certain structured settings.
We refer to this method as Structured Nuclear Norm Minimization (Structured NNM). 

Equation~\eqref{eq:structured nnm} very closely resembles the problem of decomposing a matrix into a low-rank component and a sparse component (see e.g.~\cite{chandrasekaran2011rank}).
A popular method is Robust Principal Component Analysis (RPCA)~\cite{candes2011robust}, where one assumes that a low-rank matrix has some set of its entries corrupted.
In a more recent paper~\cite{rontogiannis2020online}, reweighted least squares optimization is applied to develop a novel online
Robust PCA algorithm for sequential data processing.

\subsection{Proposed Algorithm: Structured IRLS}
\label{ss:Structured IRLS}

We propose an iterative reweighted least squares algorithm related to~\cite{fornasier2011low,mohan2012iterative} for matrix completion with structured observations.
In particular, we adjust the IRLS-$p$ algorithm proposed in~\cite{mohan2012iterative} to take into account the sparsity-based structure in the missing entries.

The $k$-th iteration of IRLS-$p$ is given by
\begin{equation*}
    X^{k} = \argmin \limits _{X} \left \{ \|(W_p^{k-1})^{1/2}X\|_F^2 \colon  \cP_{\Omega}(X) = \cP_{\Omega} (M)  \right \},
\end{equation*}
where $X^{k}\in \R ^{m \times n}$ denotes the $k$-th approximation of the true matrix $M$, $W _p ^ k \in \R^{m \times m}$ is a weight matrix defined as $W _p ^ 0 = I$, and for $k>0$, $W _p ^ k = ((X^{k})^{\top}X^k + \gamma ^k I)^{\frac{p}{2}-1}$.
Here $\gamma^k > 0$ is a regularization parameter added to ensure that $W _p ^ k$ is well-defined.

We adjust IRLS-$p$ by adding a regularization term on the unobserved entries in each iteration, namely a weighted $\ell_2$ norm as proposed in~\cite{daubechies2010iteratively} to induce sparsity.
We define the corresponding weights at the $k$-th iteration as
$w_q^{k}= (z^2(X^k)+ \epsilon^k \mathbbm 1)^{\frac{q}{2}-1},$
where $0 < \epsilon^{k} \leq \epsilon^{k-1}$, and $0\leq q \leq 1$.
Here $z(X^k)$ denotes the vector of missing entries of the the $k$-th approximation $X^k$, and recall that $z^2(X^k)$ denotes the vector with entries squared. 
The algorithm is then designed to promote low-rank structure in the recovered matrix with sparsity in the missing entries at each iteration. 
We give a description of the choice of parameters in Section~\ref{ss: choice of parameters}.

In many applications, missing values tend to be near a certain value, e.g.\ the maximum possible value in the range, or alternatively the lowest possible value (``1 star" in movie ratings).
In cases where this value is nonzero, our objective function can be adjusted accordingly.
We refer to the algorithm as \emph{Structured IRLS}; it is outlined in Algorithm~\ref{alg:structured IRLS}.
Note that each iteration of Structured IRLS solves a quadratic program, and for $\alpha = 0$, the algorithm reduces to IRLS-$p$ studied in~\cite{mohan2012iterative}.

\begin{algorithm}[htb]
\setstretch{1.5}
\SetAlgoLined
  \SetKwInOut{Input}{input}
    \SetKwInOut{Output}{output}
    \Input{$\cP_\Omega$, $M$}
        \SetKwInOut{Set}{set}
    \Set{$k=1$, $\alpha>0$, and $0 \leq p,q \leq 1$}
    \SetKwInOut{Initialize}{initialize}
    \Initialize{$X^ 0 = \cP_\Omega (M)$, $W_p^0 = I$, $w_q^0 = \mathbbm 1$, $\gamma ^ 1 > 0$, $\epsilon ^ 1 > 0$}
     \While{not converged \do}{ 
     $X^ k = \argmin \limits _{X} \left \{ \|(W_p^{k-1})^{1/2}X\|_F^2 + \alpha \|z(X)\|_{\ell_2(w_q^{k-1})}^2
     ~:~ \cP_\Omega (X)~=~\cP_\Omega (M) \right \}$ 
     
     \SetKwInOut{Compute}{compute}
     \Compute{$W _p ^ k = ((X^{k})^{\top}X^k + \gamma ^k I)^{\frac{p}{2}-1}$ and $w_q^{k}= (z^2(X^k)+ \epsilon^k \mathbbm 1)^{\frac{q}{2}-1}$}
     \SetKwInOut{Update}{update}
     \Update{$0 < \gamma^{k+1} \leq \gamma^k$, $0 < \epsilon^{k+1} \leq \epsilon^k$}
     \SetKwInOut{Set}{set}
     \Set{$k = k + 1$}

} 
\caption{Structured IRLS for Matrix Completion}
 \label{alg:structured IRLS}
\end{algorithm}

Each iteration of Structured IRLS solves a quadratic program.
The algorithm can be adjusted to have the $\ell_2$ norm for the regularization term on the unobserved entries by fixing the weights $w_q^k = \mathbbm 1$.
Further, we can impose nonnegativity constraints on the missing entries by thresholding all missing entries to be nonnegative.

We now provide an analytic remark, similar to~\cite[Proposition 1]{molitor2018matrix}, applied to the objective functions for each iteration of IRLS~\cite{mohan2012iterative} and Structured IRLS.
We consider the simplified setting in which all of the unobserved entries are exactly zero. 
We show that the approximation given by an iteration of Structured IRLS will always perform at least as well as that of IRLS with the same weights assigned.
This remark is weaker than~\cite[Proposition 1]{molitor2018matrix} as it does not apply to the entire algorithm; instead it only bounds the performance of a single iterative step.

\begin{remark}
Let \[\tilde X = \argmin \limits _{X} \left \{ \|W^{1/2}X \|_F^2 ~:~ \cP_\Omega (X) = \cP_\Omega (M) \right \}\] be the minimizer of the objective function of each iterate in IRLS~\cite{mohan2012iterative}.
Let \[\hat X = \argmin \limits _{X} \left \{ \|W^{1/2}X \|_F^2 + \alpha \|\cP_{\Omega^c}(X)\|^2 ~:~ \cP_\Omega (X) = \cP_\Omega (M) \right \}\] be the minimizer of the objective function generalizing\footnote{Here $\| \cdot \|$ is an arbitrary matrix norm; one recovers Structured IRLS by choosing the norm $\ell_2(w)$.} each iterate in Structured IRLS (with $\alpha >0$).
If $\cP_{\Omega^c} (M)$ is the zero matrix and the same weights $W$ are assigned, then $ \|M - \hat X \| \leq \|M - \tilde X \|$ for any matrix norm $\| \cdot \|$.
\end{remark}
\begin{proof}
By definition of $\hat X$, we have $\|W^{1/2}\hat X \|_F^2 + \alpha \|\cP_{\Omega^c}(\hat X)\|^2 \le \|W^{1/2}\tilde X \|_F^2 + \alpha \|\cP_{\Omega^c}(\tilde X)\|^2$.
Similarly, by definition of $\tilde X$, we have $\|W^{1/2} \tilde X \|_F^2 \leq \|W^{1/2} \hat X \|_F^2$.
Therefore, 
\begin{align*}
\|W^{1/2} \hat X \|_F^2 + \alpha \|\cP_{\Omega^c}(\hat X)\|^2 &\leq \|W^{1/2} \tilde X \|_F^2 + \alpha \|\cP_{\Omega^c}(\tilde X)\|^2 \\
&\leq \|W^{1/2} \hat X \|_F^2 + \alpha \|\cP_{\Omega^c}(\tilde X)\|^2.
\end{align*}
Since $\alpha>0$, this implies $\|\cP_{\Omega^c}(\hat X)\|^2 \leq \|\cP_{\Omega^c}(\tilde X)\|^2$.
We have
\begin{align*}
\|M - \hat X \| &= \|\cP_{\Omega^c} (\hat X) \| &&\text{since } \cP_\Omega (M) = \cP_\Omega (\hat X) \text{ and } \cP_{\Omega^c} (M)=0\\
&\le\|\cP_{\Omega^c} (\tilde X)\| \\
&=\|M - \tilde X \|  &&\text{since } \cP_\Omega (M) = \cP_\Omega (\tilde X) \text{ and } \cP_{\Omega^c} (M)=0.
\end{align*}
\end{proof}

\subsection{Proposed Implementation: Structured sIRLS}
\label{ss: Structured sIRLS}
In this section, we propose a gradient-projection-based implementation of Structured IRLS, that we will refer to as \emph{Structured sIRLS}.
Indeed, sIRLS stands for short IRLS (in analogy to~\cite{mohan2012iterative}), the reason being we do not perform gradient descent until convergence; instead we take however many steps desired.
Further, calculating $\cP_ \Omega(X)$ is computationally cheap, so the gradient projection algorithm can be used to efficiently solve the quadratic program in each iteration of Structured IRLS. 

In this implementation, we do not perform projected gradient descent on
\[ \|(W_p^{k-1})^{1/2}X\|_F^2 + \alpha \|z(X)\|_{\ell_2(w_q^{k-1})}^2, \]
with $\cP_\Omega (X) = \cP_\Omega (M)$ for each iteration $k$.
Instead, we perform projected gradient descent on $\|z(X)\|_{\ell_2(w_q^{k-1})}^2$ and $\|(W_p^{k-1})^{1/2}X\|_F^2$ consecutively. 
This allows us to update the weights before each alternating step, and to control how many gradient steps we would like to perform on each function.

We follow~\cite{mohan2012iterative} for the derivation of the gradient step of $\|(W_p^{k-1})^{1/2}X\|_F^2$ at the $k$-th iteration.
Indeed, we consider the smooth Schatten-$p$ function, for $p>0$: 
\begin{equation*}
\begin{aligned}
f_p(X) = \tr (X^\top X + \gamma I)^{\frac{p}{2}} = \sum\limits_{i=1}^n \left( \sigma_i^2(X) + \gamma \right)^{\frac{p}{2}}.
\end{aligned}
\end{equation*} 
Note that $f_p(X)$ is differentiable for $p > 0$, and convex for $p \geq 1$~\cite{mohan2012iterative}.
For $\gamma = 0$ we have $f_1(X) = \|X\|_*$,
which is also known as the Schatten-1 norm.
Again for $\gamma = 0$, we have $f_p(X) \to \rank(X)$ as $p \to 0$~\cite{mohan2012iterative}.
Further, for $p=0$, we define
\[\log\det(X^\top X + \gamma I),\]
a smooth surrogate for $\rank(X^\top X)$ (see e.g.~\cite{fazel2002matrix,fazel2003log,mohan2012iterative,recht2010guaranteed}).
Thus, it is of interest to minimize $f_p(X)$ subject to the set of constraints $\cP_\Omega (X) = \cP_\Omega (M)$ on the observed entries.

The gradient projection iterates of Structured sIRLS are given by
\begin{equation*}
X^{k+1}  = \cP_{\Omega^c}(X^k - s^k\nabla f_p(X^k)) + \cP_{\Omega}(M),
\end{equation*} 
where $s^k$ denotes the gradient step size at the $k$-th iteration and $\nabla f_p(X^k) = X^k W_p^k$,
where
we iteratively define $W_p^k$ as
\begin{equation*}
W_p^{k}  = ({X^k}^\top X^k + \gamma^k I)^{\frac{p}{2}-1},
\end{equation*} 
with $0 < \gamma^{k} \leq \gamma^{k-1}$.
This iterate describes our gradient step promoting low-rankness, where we preserve the observed entries and update only the missing entries.

Further, we promote sparsity in the missing entries as follows.
Instead of minimizing the $\ell_1$ norm of the vector of missing entries, we iteratively minimize a re-weighted $\ell_2$ norm of missing entries as described in~\cite{daubechies2010iteratively}.
Let $z(X^k)$ denote the vector of missing entries of the the $k$-th approximation $X^k$.
Define the weighted $\ell_2$ norm of $z(X)$ as
\[g_q(X) = \|z(X)\|_{\ell_2(w_q)}^2  = \sum\limits_{i=1}^{mn-|\Omega|} (w_q)_i z_i^2(X), \]
where $(w_q)_i = ( z_i^2(X) + \epsilon)^{q/2-1}$ (as done in~\cite{daubechies2010iteratively}).
The $i$-th entry of the gradient of $g_q(X)$ is given by $(\nabla g_q(X))_i = 2(w_q)_i z_i$.
Therefore, the gradient projection iterates are given by
\begin{equation*}
    z(X^{k+1}) = z(X^{k}) -c^k\nabla g_q(X^k),
\end{equation*}
where $c^k$ denotes the gradient step size at the $k$-th iteration.
We iteratively define the weights $w_q^k$ as
\begin{equation*}
w_q^{k}= (z^2(X^k)+ \epsilon^k \mathbbm 1)^{\frac{q}{2}-1},
\end{equation*} 
where $0 < \epsilon^{k} \leq \epsilon^{k-1}$.

We outline in Algorithm~\ref{alg:structured sIRLS} Structured sIRLS, a gradient-projection-based implementation of Structured IRLS.

\begin{algorithm}[H]
\setstretch{1.35}
\SetAlgoLined
  \SetKwInOut{Input}{input}
    \Input{$\cP_\Omega$, $M$, $r$}
        \SetKwInOut{Set}{set}
    \Set{$k=1$, $0 \leq p,q \leq 1$, $k_s>0$, $k_l >0$, $c^k>0$, $s^k > 0$}
    \SetKwInOut{Initialize}{initialize}
    \Initialize{$X^ 0 = \cP_\Omega (M)$, $w_q^0 = \mathbbm 1$, $\gamma ^ 1 > 0$, $\epsilon ^ 1 > 0$}
     \While{not converged}{
     \SetKwInOut{Find}{find}
     \SetKwInOut{Compute}{compute}
     \SetKwInOut{Perform}{perform}
     \SetKwInOut{Update}{update}
    \Perform{take $k_s$ steps promoting sparsity, $z(X^{k}) = z(X^{k-1}) -c^k({w_q}^{k-1} \odot {z(X^{k-1})})$}
    \Update{update the weights promoting low-rankness, $W_p^{k}  = ({X^k}^\top X^k + \gamma^k I)^{\frac{p}{2}-1}$}
    \Perform{take $k_l$ steps promoting low-rankness, $X^{k+1}  = P_{\Omega^c}(X^k - s^kX^k W_p^k) + P_{\Omega}(M)$}
    \Update{update the weights promoting sparsity, $w_q^{k}= (z^2(X^{k+1})+ \epsilon^k \mathbbm 1)^{\frac{q}{2}-1}$
    }
     \Update{update the regularizers, $0 < \gamma^{k+1} \leq \gamma^k$, $0 < \epsilon^{k+1} \leq \epsilon^k$}
    \Set{set $k = k + 1$}
}
 \caption{Structured sIRLS for Matrix Completion}
 \label{alg:structured sIRLS}
\end{algorithm}

A rank estimate $r$ of the matrix $M$ is used as an input to truncate the singular value decomposition (SVD) when computing the weights $W_p^k$.
In our implementation, we use a randomized algorithm for SVD computation~\cite{halko2011finding}.
When the rank of the matrix is not estimated or provided, we instead choose $r$ to be $\min\{r_{max}, \hat r\}$ at each iteration, where $\hat r$ is the largest integer such that $\sigma_{ \hat r}(X^k ) > 10^{-2} \cdot \sigma_1(X^k )$, and where $r_{max} = \left \lceil n \left ( 1 - \sqrt{1 - \frac{| \Omega|}{mn}}\right) \right \rceil$ (as implemented in~\cite{mohan2012iterative}).

\section{Numerical Experiments}
\label{sec: numerical experiments}

In this section, we run numerical experiments to evaluate the performance of Structured sIRLS.
We compare Structured sIRLS to the performance of sIRLS (studied in~\cite{mohan2012iterative}) and Structured NNM (studied in~\cite{molitor2018matrix}) on structured settings.
Our code for Structured sIRLS is available at~\cite{StructuredIRLScode}.
Further, we use the publicly available code of sIRLS~\cite{IRLScode}.

First, in Section~\ref{ss: choice of parameters}, we explain the choice of parameters we use.
We describe our experiments for exact matrix completion in Section~\ref{sec: exact recovery}.
For problems of varying degrees of difficulty in terms of the sampling rate, degrees of freedom, and sparsity levels, we find that Structured sIRLS often outperforms sIRLS and Structured NNM in the structured setting.
In Section~\ref{sec:noisy completion} we consider matrix completion with noise, finding that Structured sIRLS improves upon sIRLS in the structured setting with low noise.
As the noise level increases, the performance of Structured sIRLS remains controlled, approximately the same as the performance of sIRLS.

\subsection{Choice of parameters} 
\label{ss: choice of parameters}
In all the numerical experiments, we adopt the same parameters. 
However, one can use different choices for parameters, or optimize some of the parameters.
We normalize the input data matrix $M$ to have a spectral norm of 1 (as done in~\cite{mohan2012iterative}).

We are particularly interested in the case $p=q=1$. In our experiments, we set $p = q = 1$, but generally, these parameters can be varied over the range $0 \leq p,q \leq 1$.
Each value of $p$ and $q$ define a different objective function (see Sections~\ref{ss:related sparse recovery} and~\ref{ss:related low-rank matrix recovery}).

For the implementation parameters, we set $k_s = 1$ and $k_l =10$, which means that we take one gradient step to promote sparsity and ten gradient steps to promote low-rankness, respectively. 
These parameters can be varied based on the low-rankness of the matrix and on  the expected sparsity of its missing entries.
Further, we set the regularizers $\gamma^k = (1/2)^k$ and $\epsilon^k = (9/10)^k$ at the $k$-th iteration. 
However, there are other possible choices for these regularizers, for example $\epsilon^k$ could depend on the $(s+1)$-th largest value of $z(X^k)$, where $s$ is the sparsity of $z(X^k)$ (as done in~\cite{daubechies2010iteratively}).
Similarly, $\gamma^k$ could depend on the $(r+1)$-th singular value of $X^k$, where $r$ is the rank of $M$ (as done in~\cite{fornasier2011low}).

Lastly, for all $k$ we set the step size $s^k = (\gamma ^k)^{1-\frac{p}{2}}$ to promote low-rankness and $c^k = 10^{-6}$ to promote sparsity; however, these parameters could be scaled or varied.
We define the relative distance between two consecutive approximations as
\[ d(X^k, X^{k-1}) = \|X^k - X^{k-1} \|_F / \|X^k\|_F. \]
We say the algorithm converges if we obtain $d(X^k, X^{k-1}) < 10^{-5}$.
We set the tolerance $10^{-5}$ for both sIRLS and Structured sIRLS in our comparison experiments,\footnote{In the original implementation of sIRLS provided by the authors~\cite{IRLScode,mohan2012iterative}, the tolerance value is set to $10^{-3}$. 
However, Structured sIRLS converges much faster per iteration, thus attaining the tolerance $10^{-3}$ with fewer iterations.
To report fair comparisons between the algorithms that do not overly benefit Structured sIRLS, we set the tolerance to $10^{-5}$ in addition to increasing the maximum number of iterations for sIRLS.}
and we set the maximum number of iterations for Structured sIRLS to be 1000 and for sIRLS to be 5000.

\subsection{Exact Matrix Completion}
\label{sec: exact recovery}
We first investigate the performance of the Structured sIRLS algorithm when the observed entries are exact, i.e.\ there is no noise in the observed values.
We construct $m\times n$ matrices of rank $r$ as done in~\cite{molitor2018matrix}.
We consider $M = M_L M_R$,
where $M_L \in \R^{m\times r}$ and $M_R \in \R^{r\times n}$ are sparse matrices.
Indeed, the entries of $M_L$ (resp.\ $M_R$) are chosen to be zero uniformly at random so that on average $70\%$ (resp.\ $50\%$) of its entries are zero.
The remaining nonzero entries are uniformly distributed at random between zero and one.
The sparsity level of the resulting matrix $M$ cannot be calculated exactly from the given sparsity levels of $M_L$ and $M_R$.
Thus for each of the following numerical simulations, we indicate on average the sparsity level of $M$ (we refer to the density of $M$ as the fraction of nonzero entries).

\begin{figure}[h]
\centering
\includegraphics[width=4.3in]{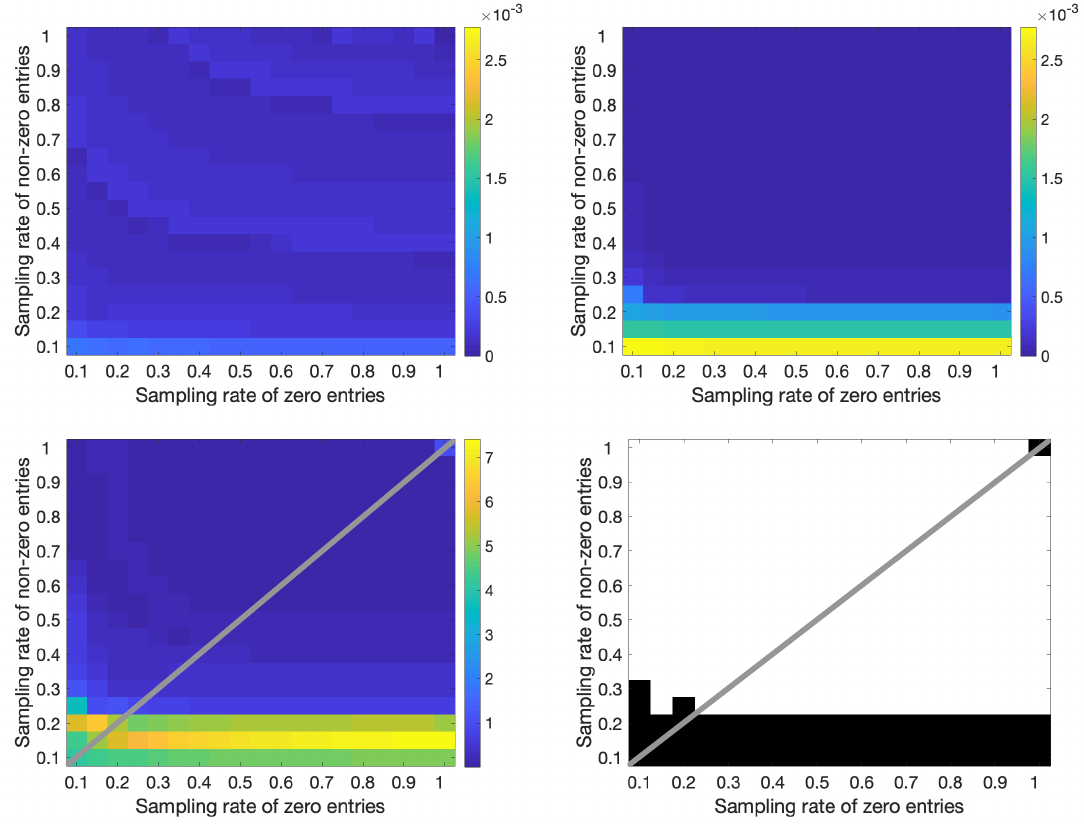}
\caption{We consider twenty $1000 \times 1000$ sparse random matrices of rank 10
with average density equal to 0.66.
The upper plots display
(left) the average relative error for sIRLS $\|M - \tilde X\|_F/\| M\|_F$, and (right) the average relative error for Structured sIRLS $\|M - \hat X\|_F/\| M\|_F$.
The lower plots display
(left) the average ratio $\|M - \hat X\|_F / \|M - \tilde X\|_F$, and
(right) the binned average ratio where white means the average ratio is strictly less than 1, and black otherwise.
}
\label{fig:exp1000_rank10}
\end{figure}

\begin{figure}[h]
\centering
\includegraphics[width=4.3in]{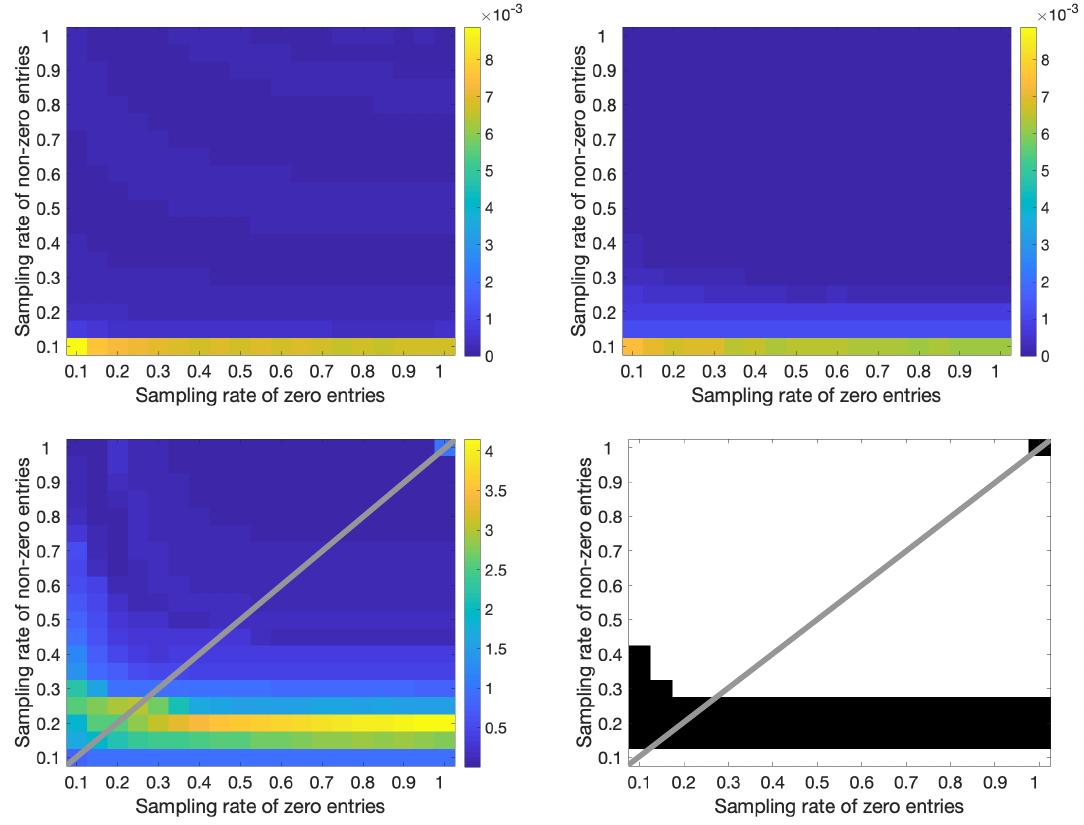}
\caption{We consider twenty $500 \times 500$ sparse random matrices of rank 10 with average density equal to 0.66.
The upper plots display
(left) the average relative error for sIRLS $\|M - \tilde X\|_F/\| M\|_F$, and (right) the average relative error for Structured sIRLS $\|M - \hat X\|_F/\| M\|_F$.
The lower plots display
(left) the average ratio $\|M - \hat X\|_F / \|M - \tilde X\|_F$, and
(right) the binned average ratio where white means the average ratio is strictly less than 1, and black otherwise.
}
\label{fig:exp500_rank10}
\end{figure}

For each experiment with $m$, $n$, and $r$ fixed, we choose twenty random matrices of the form $M=M_L M_R$.
We subsample from the zero and nonzero entries of the data matrix at various rates to generate a matrix with missing entries.
We define the \emph{relative error of Structured sIRLS} as
\[\|M - \hat X\|_F/\| M\|_F,\]
where $\hat X$ is the output of the Structured sIRLS algorithm. 
Similarly, we define the \emph{relative error of  sIRLS} as
\[\|M - \tilde X\|_F/\| M\|_F,\]
where $\tilde X$ is the output of the sIRLS algorithm. 
The \emph{average ratio} is then defined as
\[\|M - \hat X\|_F / \|M - \tilde X\|_F.\]
We say Structured sIRLS outperforms sIRLS when the average ratio is less than one, and vice versa when the average ratio is greater than or equal to one.
These two cases, when the average ratio is strictly less than or greater than or equal to one, are visually represented by the white and black squares, respectively, in the bottom right plots of Figures~\ref{fig:exp1000_rank10}--\ref{fig:exp100_rank10} and~\ref{fig:exp100_rank20}.
We refer to this binary notion of average ratio as \textit{binned average ratio}.
We report each of these error values in our numerical experiments.

It is important to note that the setting we are interested in is the structured setting where the submatrix of missing values is close to 0 in the $L_1$ or $L_0$ norm sense. 
This setting can be observed in the upper left triangle of the images
in Figures~\ref{fig:exp1000_rank10}--\ref{fig:exp100_rank3wnoise4}
(in particular, this is the region above the diagonal gray line in the bottom rows of Figures~\ref{fig:exp1000_rank10}--\ref{fig:exp100_rank3wnoise4}).
In this upper-left triangular region, the percentage of nonzero entries that are sampled is greater than the percentage of zero entries that are sampled.
Hence the region above the diagonal gray lines is the structured setting that Structured sIRLS is designed for.

In general, algorithms obtain better accuracy as we move right along a row or up along a column in Figures~\ref{fig:exp1000_rank10}--\ref{fig:exp100_rank3wnoise4}, since we are sampling more and more entries.
In addition, it is important to note that in all experiments we are using the same algorithm (with fixed parameters) for all the cases considered in our computations, without any parameter optimization.
The Structured sIRLS algorithm promotes sparsity in all the cases, even in the unstructured settings.
Omitting the sparsity promoting step would result in an algorithm promoting only low-rankness.

\subsubsection{$1000 \times 1000$ rank 10 matrices}
In Figure~\ref{fig:exp1000_rank10}, we construct twenty random matrices of size $1000 \times 1000$ and of rank 10, as described in Section~\ref{sec: exact recovery}.
Error ratios below one in the bottom left plot of Figure~\ref{fig:exp1000_rank10} indicate that Structured sIRLS outperforms sIRLS. 
In this particular experiment, we observe that Structured sIRLS outperforms sIRLS for most of the structured cases (the upper left triangle above the gray line), and more.
For this particular experiment, it turns out that this happens roughly when the decimal percentage of sampled nonzero entries is greater than 0.2.

Note that in the case where all entries are observed (no longer a matrix completion problem), both relative errors are 0 and thus the average ratio is 1.
We only say that Structured sIRLS outperforms sIRLS when the average ratio is strictly less than 1, and this is why the upper right pixel in the bottom right plot of Figure~\ref{fig:exp1000_rank10} is black.
The same is true in later figures.


\subsubsection{$500 \times 500$ rank 10 matrices}
In Figure~\ref{fig:exp500_rank10}, we construct twenty sparse random matrices of size $500 \times 500$ and of rank 10, as described in Section~\ref{sec: exact recovery}.
We observe that Structured sIRLS outperforms sIRLS not only in the majority of the structured cases, but also in many of the other cases where the submatrix of unobserved entries is not necessarily sparse.

\subsubsection{$100 \times 100$ rank 10 matrices}
\begin{figure}[h]
\centering
\includegraphics[width=4.3in]{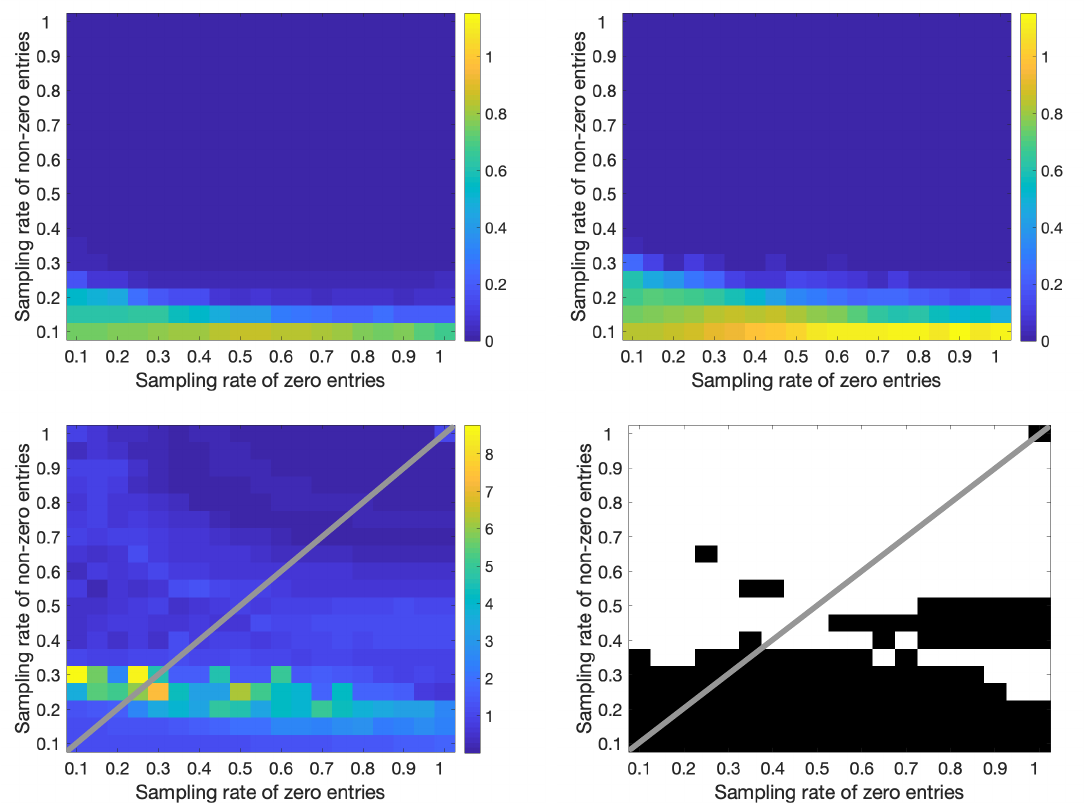}
\caption{We consider twenty $100 \times 100$ sparse random matrices of rank 10 with average density equal to 0.66.
The upper plots display
(left) the average relative error for sIRLS $\|M - \tilde X\|_F/\| M\|_F$, and (right) the average relative error for Structured sIRLS $\|M - \hat X\|_F/\| M\|_F$.
The lower plots display
(left) the average ratio $\|M - \hat X\|_F / \|M - \tilde X\|_F$, and
(right) the binned average ratio where white means the average ratio is strictly less than 1, and black otherwise.
}
\label{fig:exp100_rank10}
\end{figure}

In Figure~\ref{fig:exp100_rank10}, we construct twenty random matrices of size $100 \times 100$ and of rank 10, as described in Section~\ref{sec: exact recovery}.
We observe in Figure~\ref{fig:exp100_rank10} that Structured sIRLS outperforms sIRLS when 
the sampling rate of the nonzero entries is high (roughly speaking, when the decimal percentage of sampled nonzero entries is greater than 0.5), which covers the majority of the cases where there is sparsity-based structure in the missing entries.

\subsubsection{$100 \times 100$ matrices with no knowledge of the rank a priori}
\begin{figure}[h]
\centering
\includegraphics[width=4.3in]{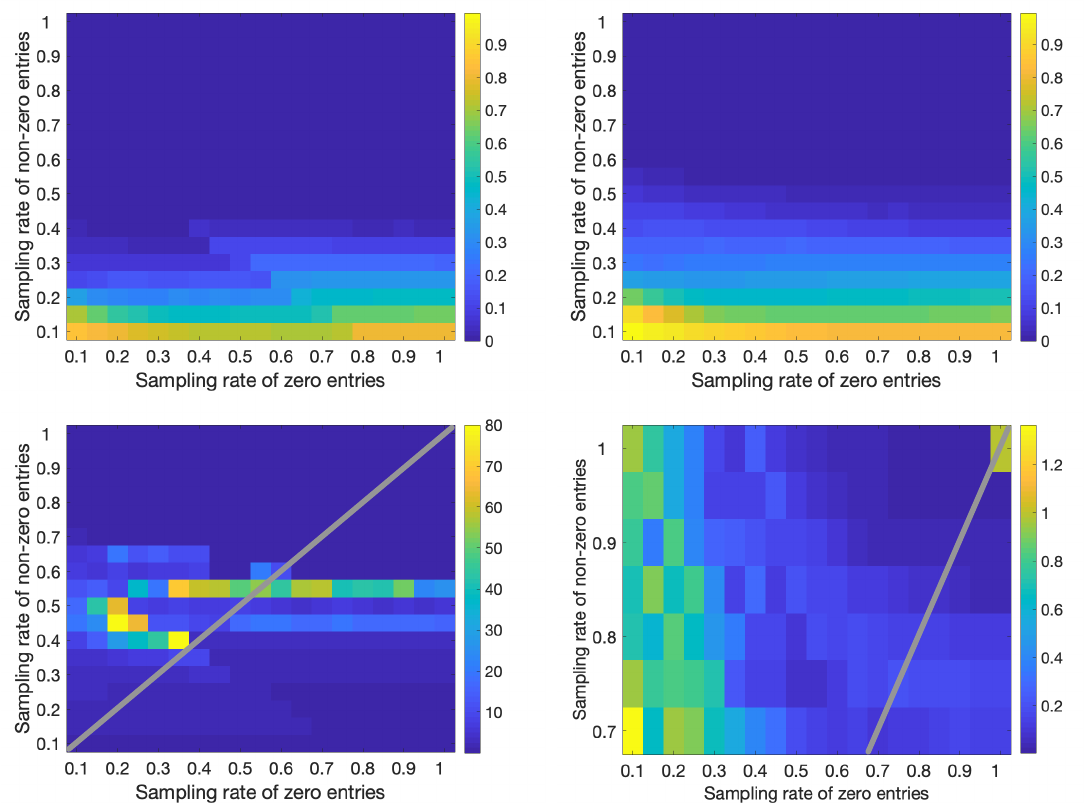}
\caption{We consider twenty $100 \times 100$ sparse random matrices of rank 8 with density equal to 0.58, but we do not input any rank guess.
The upper plots display
(left) the average relative error for sIRLS $\|M - \tilde X\|_F/\| M\|_F$, and (right) the average relative error for Structured sIRLS $\|M - \hat X\|_F/\| M\|_F$.
The lower plots display
(left) the average ratio $\|M - \hat X\|_F / \|M - \tilde X\|_F$,
and 
(right) the average ratio when the sampling rate of non-zero entries is at least 0.70 (a zoomed in version of part of the lower left plot).
}
\label{fig:exp100_rank8tounkown}
\end{figure}

In Figure~\ref{fig:exp100_rank8tounkown}, we construct twenty random matrices of size $100 \times 100$ and of rank 8, as described in Section~\ref{ss: Structured sIRLS}.
For this experiment, we do not provide the algorithm with any rank estimate, for either sIRLS or Structured sIRLS.
Instead, we allow the algorithm to estimate the rank at each iteration based on a heuristic described in Section~\ref{ss: Structured sIRLS}.
We observe in the bottom right plot of Figure~\ref{fig:exp100_rank8tounkown}, 
where we zoom in on the cases where the sampling rate of non-zero entries is at least 0.7, that Structured sIRLS outperform sIRLS to some extent in this region.
Indeed, Structured sIRLS does particularly better when more entries are observed.

\subsubsection{$100 \times 100$ rank 20 matrices}
\label{exp100_rank20}

We say a matrix completion problem is hard when the degrees of freedom ratio $FR$ is greater than $0.4$ (as in~\cite{mohan2012iterative}).
In the previous experiments, we considered a few cases where $FR > 0.4$, which occur when the sampling rates of zero and nonzero entries are both relatively small.
In these cases, there is not necessarily high sparsity-based structure, which imposes another challenge since the sampling rate of non-zero entries is approximately equal to or only slightly greater than the sampling rate of zero entries.
Therefore, in this section, we consider hard cases (where $FR >0.4$) with sparsity-based structure.

\begin{figure}[h]
\centering
\includegraphics[scale=0.14]{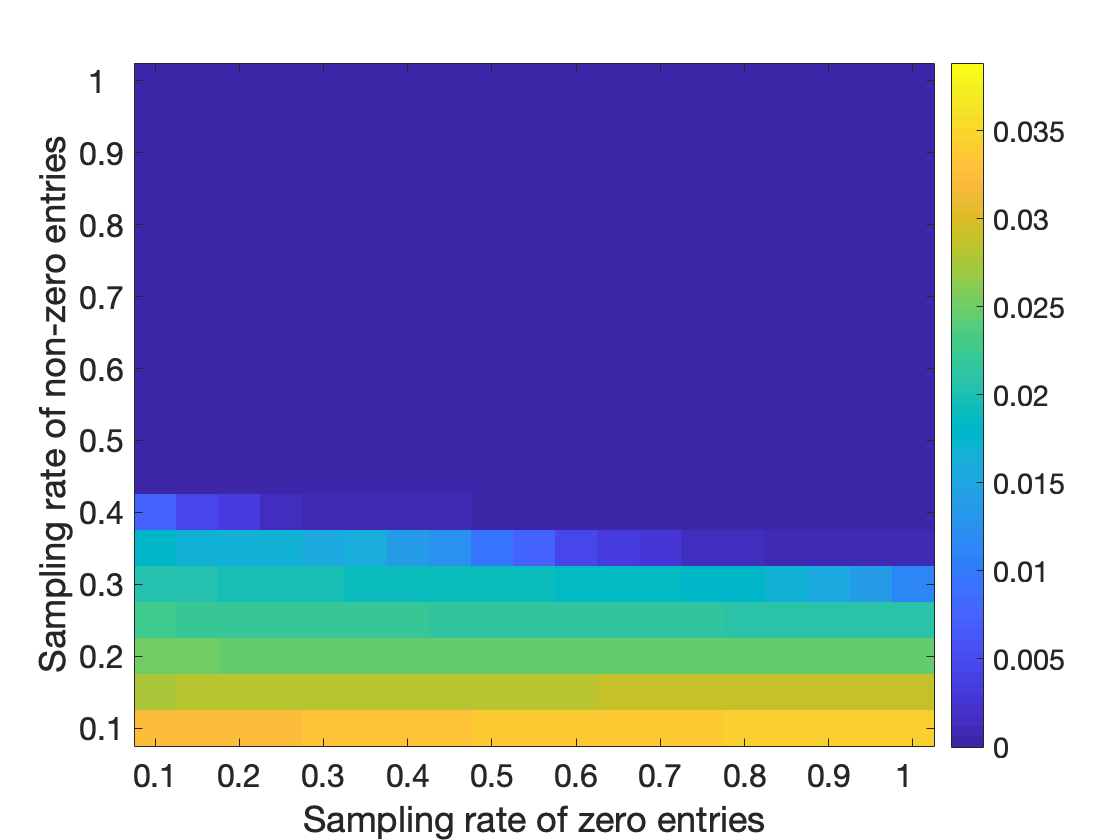}
\includegraphics[scale=0.14]{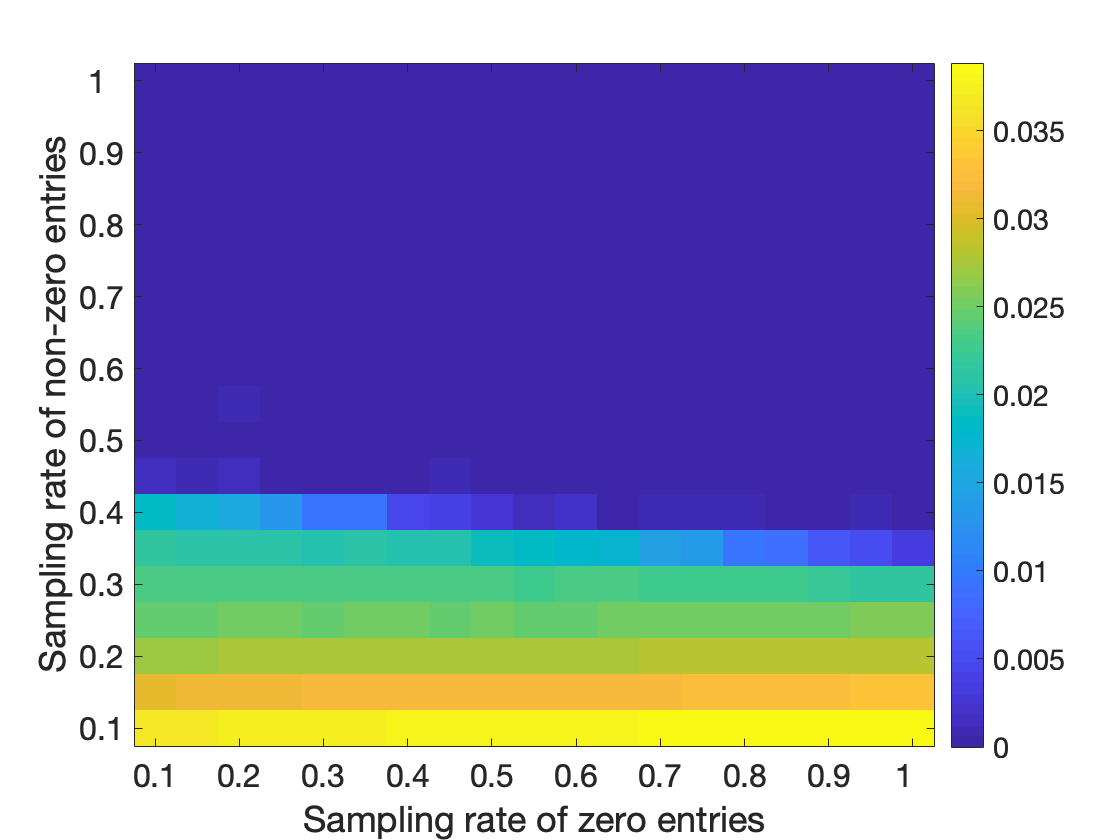}

\includegraphics[width=2.2in]{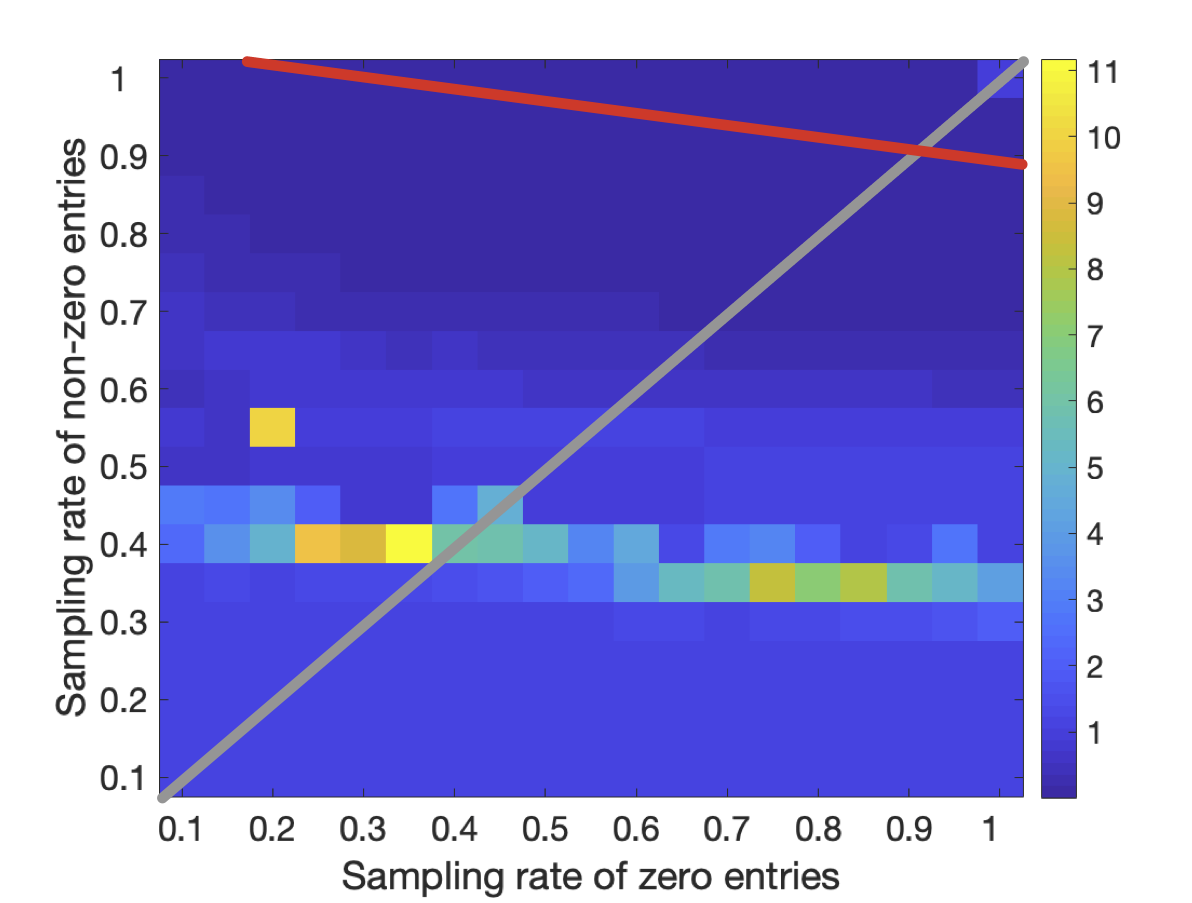}
\includegraphics[width=2.2in]{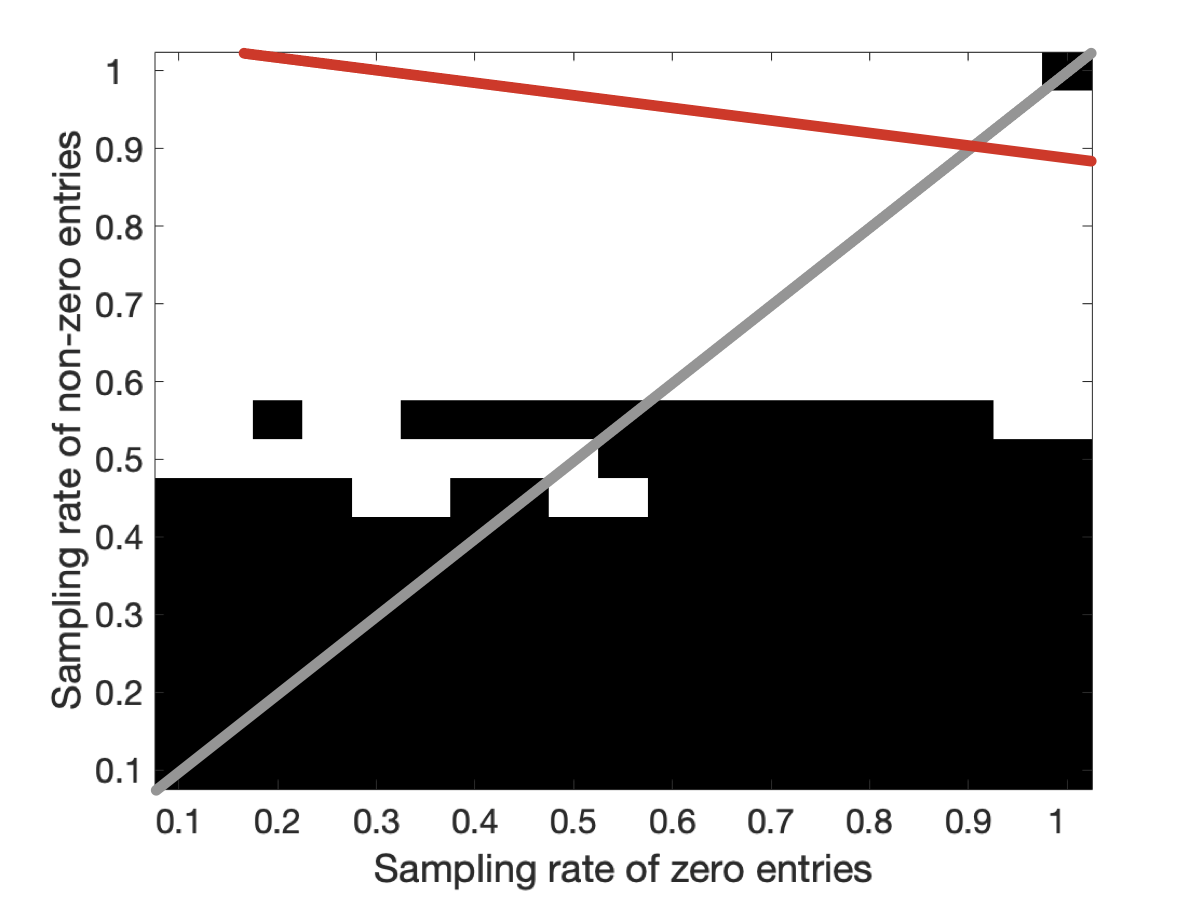}
\caption{We consider twenty $100 \times 100$ sparse random matrices of rank 20 with average density equal to 0.88.
The upper plots display
(left) the average relative error for sIRLS $\|M - \tilde X\|_F/\| M\|_F$, and (right) the average relative error for Structured sIRLS $\|M - \hat X\|_F/\| M\|_F$.
The lower plots display
(left) the average ratio $\|M - \hat X\|_F / \|M - \tilde X\|_F$, and
(right) the binned average ratio where white means the average ratio is strictly less than 1, and black otherwise. 
The red line separates the hard cases from those that are not: 
all the cases below the red line are hard.
}
\label{fig:exp100_rank20}
\end{figure}

In Figure~\ref{fig:exp100_rank20}, we construct twenty random matrices of size $100 \times 100$ and of rank 20, as described in Section~\ref{sec: exact recovery}.
If the number of sampled entries is 90\% of the entire matrix, i.e.\ $|\Omega| = 0.9 \cdot m \cdot n$, then 
\[FR = r(m+n-r) / |\Omega| = 20(200-20) / (0.9 \times 100^2) =  0.4.\]
So, even sampling 90\% of the matrix is still considered to be a hard problem.
In the bottom row of Figure~\ref{fig:exp100_rank20}, the added red line separates the ``hard" cases from those that are not: 
all the cases below the red line are hard.
Note that in these hard regimes with sparsity-based structure, Structured sIRLS outperform sIRLS more often than not.

\subsubsection{Comparison with Structured NNM on $30\times 30$ rank 7 matrices}
\label{ss: comparison to Struc. NNM}
\begin{figure}[h!]
\centering
\includegraphics[width=4.3in]{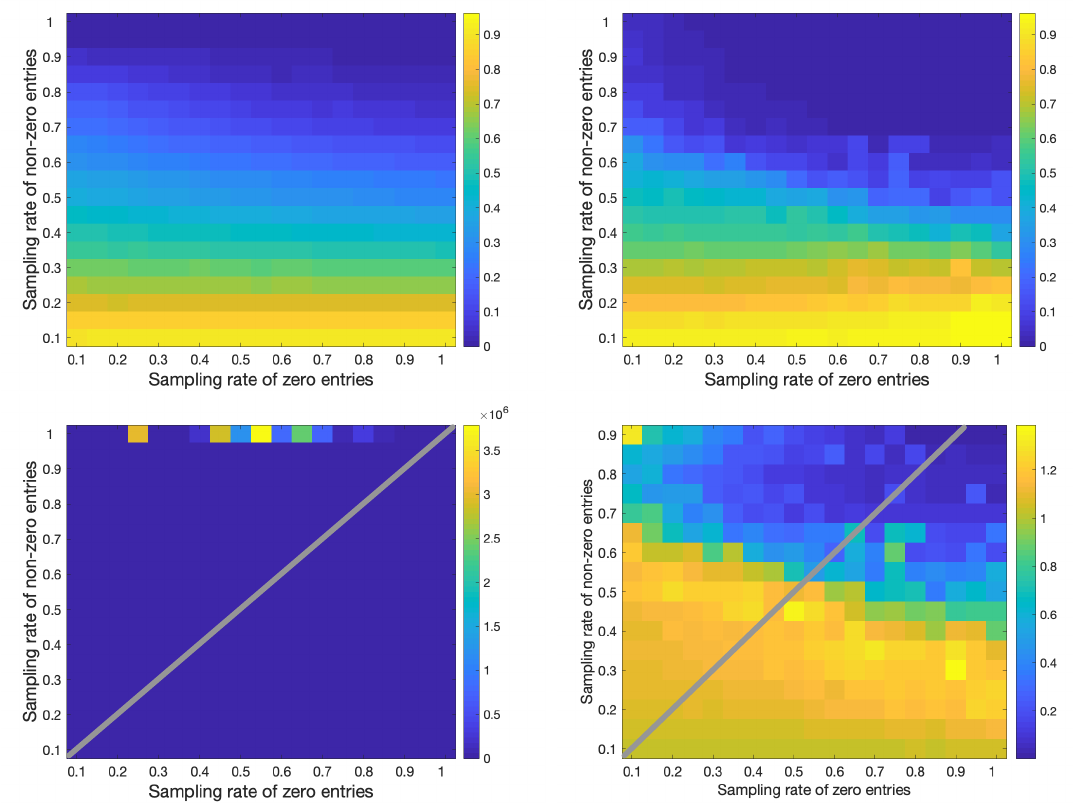}
\caption{We consider twenty $30 \times 30$ sparse random matrices of rank 7, with average density equal to 0.53.
We provide Structured sIRLS with the rank of the matrices.
We optimize for each matrix and combination of sampling rates the regularization parameter $\alpha \in \{ 10^{-4}, 10^{-3}, 10^{-2}, 10^{-1}\}$ for Structured NNM.
The upper plots display
(left) the average relative error for Structured NNM $\|M - \bar X\|_F/\| M\|_F$, and (right) the average relative error for Structured sIRLS $\|M - \hat X\|_F/\| M\|_F$.
The lower plots display
(left) the average ratio $\|M - \hat X\|_F / \|M - \bar X\|_F$,
and 
(right) the average ratio when the sampling rate of non-zero entries is at most 0.90 (a zoomed in version of part of the lower left plot).
}
\label{fig:exp30_rank7tokown_opt}
\end{figure}

In this section, we run numerical experiments to compare the performance of Structured sIRLS with Structured NNM, using the $L_1$ norm on the submatrix of unobserved entries.
We use CVX, a package for specifying and solving convex programs~\cite{gb08,cvx}, to solve Structured NNM.
In the experiments of Figure~\ref{fig:exp30_rank7tokown_opt}, we construct twenty random matrices of size $30\times 30$ and of rank 7 as described in Section~\ref{sec: exact recovery}.
We compare the accuracy with Structured NNM on small-sized matrices due to computational constraints of CVX: with this implementation of Structured NNM, it is difficult to handle significantly larger matrices.
Similar experiments are considered in~\cite{molitor2018matrix}, where Structured NNM is compared to NNM.
Comparing our iterative algorithm to Structured NNM is important since Structured NNM adapts nuclear norm minimization and $\ell_1$ norm minimization, which are common heuristics for minimizing rank and inducing sparsity, respectively.
We define the \emph{relative error of Structured NNM} as
$\|M - \bar X\|_F/\| M\|_F,$
where $\bar X$ is the output of the Structured NNM algorithm. 
The average ratio is then defined as
$\|M - \hat X\|_F / \|M - \bar X\|_F$, where $\hat X$ is the output of the Structured sIRLS algorithm.

For all sampling rates, the degrees of freedom ratio is greater than 0.4, i.e.\ all the problems are considered to be ``hard" matrix completion problems.
In Figure~\ref{fig:exp30_rank7tokown_opt}, we provide Structured sIRLS with the rank of the matrices. 

In Figure~\ref{fig:exp30_rank7tokown_opt}, we give Structured NNM an advantage by optimizing for each matrix and combination of sampling rates the regularization parameter $\alpha \in \{ 10^{-4}, 10^{-3}, 10^{-2}, 10^{-1}\}$ for Structured NNM.
However, for Structured sIRLS (again with $p=q=1$) we do not optimize the gradient step sizes or the number of step sizes. 
Varying the number or sizes of the gradient steps controls how much we would like to promote low-rankness versus sparsity in the submatrix of missing entries.
In the experiments of Figure~\ref{fig:exp30_rank7tokown_opt}, we observe that for the most part where the sampling rate of nonzero entries is between 0.6 and 0.9, Structured sIRLS performs better than Structured NNM.
Furthermore, for the remainder of the structured settings, Structured sIRLS performs approximately the same as Structured NNM or only slightly worse. 
We note that in a couple of cases where the sampling rate of nonzero entries is 1, and where the relative error for both algorithms is close to zero, Structured NNM performs much better.
This is in part because we optimize Structured NNM over $\alpha \in \{ 10^{-4}, 10^{-3}, 10^{-2}, 10^{-1}\}$;
see~\cite{molitor2018matrix} where the relationship between the choice of $\alpha$ and the accuracy of Structured NNM is studied.

As observed in Figure~\ref{fig:exp30_rank7tokown_opt}, our proposed method is often comparable to Structured NNM on small-sized matrices, with certain regions where Structured sIRLS results in improved accuracy.
In addition, iterative methods are well-known to offer ease of implementation and reduced computational resources, making our approach attractive not only in the setting of small-sized matrices, but also large-sized matrices.


\subsection{Matrix Completion with Noise}
\label{sec:noisy completion}
In this section, we investigate the performance of Structured sIRLS when the observed entries are corrupted with noise.
In particular, we compare the performance of Structured sIRLS with the performance of sIRLS.
We adapt sIRLS and Structured sIRLS for noisy matrix completion by replacing the observed entries $P_\Omega(M)$ with the noisily observed entries $P_\Omega(B)$ in the constraints, where $M$ is an unknown low-rank matrix that we wish to recover, where $P_\Omega(Z)$ is the measurement noise, and where the noisy matrix $B$ satisfies $P_\Omega(B) = P_\Omega(M) +P_\Omega(Z)$.
The algorithms for matrix recovery do not update the noisily observed entries, only the missing entries.
We define our noise model such that $\| P_\Omega(Z)\| _F = \epsilon \| P_\Omega(M)\| _F$ for a noise parameter $\epsilon$.
We do so by adding noise of the form
\[Z_{ij} = \epsilon \cdot\frac{ \| P_\Omega (M)\|_F}{\| P_\Omega (N)\|_F} \cdot N_{ij},\]
where $N_{ij}$ are i.i.d.\ Gaussian random variables with the standard distribution $\mathcal N(0, 1)$.
We define the relative error of Structured sIRLS as
$ \|B - \hat X\|_F/\| B\|_F,$
where $\hat X$ is the output of the Structured sIRLS algorithm.
Similarly, we define the relative error of sIRLS as
$\|B - \tilde X\|_F/\| B\|_F,$
where $\tilde X$ is the output of the sIRLS algorithm.

\begin{figure}[h]
\centering
\includegraphics[width=4.3in]{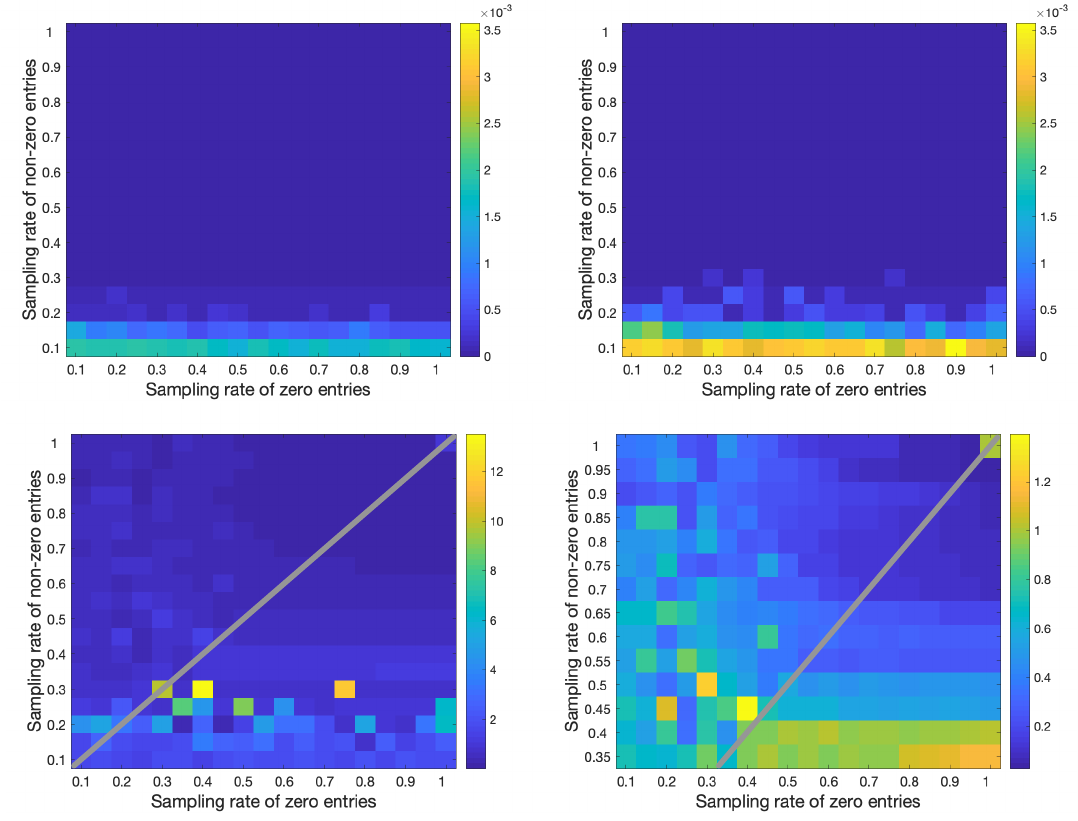}
\caption{We consider twenty $100 \times 100$ random matrices of rank 3 with noise parameter $\epsilon = 10^{-4}$.
The upper plots display
(left) the average relative error for sIRLS $\|B - \tilde X\|_F/\| B\|_F$, and (right) the average relative error for Structured sIRLS $\|B - \hat X\|_F/\| B\|_F$.
The lower plots display (left) the average ratio $\|B - \hat X\|_F / \|B - \tilde X\|_F$, and (right) the average ratio when the sampling rate of non-zero entries is at least 0.35 (a zoomed in version of part of the lower left plot).
} 
\label{fig:exp100_rank3wnoise4}
\end{figure}



In Figure~\ref{fig:exp100_rank3wnoise4}, we consider twenty random $100 \times 100$ rank 3 matrices with noise parameter $\epsilon = 10^{-4}$, where we construct our matrices in the same fashion as in Section~\ref{sec: exact recovery}.
We consider analogous structured settings as in the prior experiments, and observe that for the cases where the sampling rate of nonzero entries is greater than 0.3, which covers the  majority of the  cases  where  there  is sparsity-based structure in the missing entries, Structured sIRLS performs better than sIRLS.

For a higher noise level $\epsilon = 10^{-3}$, we observe that sIRLS and Structured sIRLS algorithms perform roughly the same.
This suggest that both sIRLS and Structured sIRLS are robust to noise, with the improvements of Structured sIRLS from the structure diminishing as the noise grows.


\section{Conclusion}

In this paper, we consider the notion of structured matrix completion, studied in the recent paper~\cite{molitor2018matrix}.
In particular, we are interested in sparsity-based structure in the missing entries whereby the vector of missing entries is close in the $\ell_0$ or $\ell_1$ norm sense to the zero vector (or more generally, to a constant vector).
For example, a missing rating of a movie might indicate the user's lack of interest in that movie, thus suggesting a lower rating than otherwise expected.
In~\cite{molitor2018matrix}, Molitor and Needell propose adjusting the standard nuclear norm minimization problem by regularizing the values of the unobserved entries to take into account the structural differences between observed and unobserved entries.

To our knowledge, we develop the first iterative algorithm that addresses the structured low-rank matrix completion problem, for which Structured NNM has been proposed.
We adapt an iterative algorithm,
called Structured IRLS, by adjusting the IRLS algorithm proposed in~\cite{mohan2012iterative}.
We also present a gradient-projection-based implementation, called Structured sIRLS that can handle large-scale matrices.
The algorithms are designed to promote low-rank structure in the recovered matrix with sparsity in the missing entries.

We perform numerical experiments on various structured settings 
to test the performance Structured sIRLS compared to sIRLS and Structured NNM.
We consider problems of various degrees of freedom and rank knowledge.
To generate matrices with sparsity-based structure in the missing entries, we subsample from the zero and nonzero entries of a sparse data matrix at various rates.
Indeed, we are interested in the structured cases, when the sampling rate of the zero entries is lower than the sampling rate of the nonzero entries.

Our numerical experiments show that Structured sIRLS often gives better recovery results than  sIRLS in structured settings.
Further,  for  small  enough noise our proposed method  often  performs  better  than  sIRLS  in  structured settings, and as noise gets larger both converge to the same performance.
Further, our numerical experiments show that Structured sIRLS is comparable to Structured NNM on small-sized matrices, with Structured sIRLS performing better in various structured regimes.

In future work, we hope to extend the theoretical results for Structured IRLS to more general settings. 
In the simplified setting, in which all of the unobserved entries are exactly zero, we show that the approximation given by an iteration of Structured IRLS will always perform at least as well as that of IRLS with the same weights assigned.
However, we empirically observe the stronger result that Structured sIRLS  often outperforms sIRLS in structured settings (in which algorithms are run until convergence, and in which not all missing entries are zero).
Another extension is to explore Structured IRLS for different values of $p$ and $q$, both empirically and theoretically.
Furthermore, a possible direction for future work is to extend sparsity-based structure in the missing entries to a more general notion of structure, whereby the probability that an entry is observed or not may depend on more than just the value of that entry.
For example, one could imagine that columns in a matrix corresponding to popular movies would have many entries (user ratings) filled in.
In this context, an entry might be more likely to be observed if many entries in its same column are also observed.

\bibliographystyle{plain}
\bibliography{MatrixCompletionForStructuredObservations.bib}

\end{document}